\documentclass[twoside,leqno,10pt]{amsart}
\usepackage{amsfonts}
\usepackage{amsmath}
\usepackage{amscd}
\usepackage{amssymb}
\usepackage{amsthm}
\usepackage{amsrefs}
\usepackage{latexsym}
\usepackage{mathrsfs}
\usepackage{bbm}
\usepackage{enumerate}
\usepackage{color}
\begin{document}

\newtheorem{theorem}[subsection]{Theorem}
\newtheorem{proposition}[subsection]{Proposition}
\newtheorem{lemma}[subsection]{Lemma}
\newtheorem{corollary}[subsection]{Corollary}
\newtheorem{conjecture}[subsection]{Conjecture}
\newtheorem{prop}[subsection]{Proposition}
\numberwithin{equation}{section}
\newcommand{\mr}{\ensuremath{\mathbb R}}
\newcommand{\mc}{\ensuremath{\mathbb C}}
\newcommand{\dif}{\mathrm{d}}
\newcommand{\intz}{\mathbb{Z}}
\newcommand{\ratq}{\mathbb{Q}}
\newcommand{\natn}{\mathbb{N}}
\newcommand{\comc}{\mathbb{C}}
\newcommand{\rear}{\mathbb{R}}
\newcommand{\prip}{\mathbb{P}}
\newcommand{\uph}{\mathbb{H}}
\newcommand{\fief}{\mathbb{F}}
\newcommand{\majorarc}{\mathfrak{M}}
\newcommand{\minorarc}{\mathfrak{m}}
\newcommand{\sings}{\mathfrak{S}}
\newcommand{\fA}{\ensuremath{\mathfrak A}}
\newcommand{\mn}{\ensuremath{\mathbb N}}
\newcommand{\mq}{\ensuremath{\mathbb Q}}
\newcommand{\half}{\tfrac{1}{2}}
\newcommand{\f}{f\times \chi}
\newcommand{\summ}{\mathop{{\sum}^{\star}}}
\newcommand{\chiq}{\chi \bmod q}
\newcommand{\chidb}{\chi \bmod db}
\newcommand{\chid}{\chi \bmod d}
\newcommand{\sym}{\text{sym}^2}
\newcommand{\hhalf}{\tfrac{1}{2}}
\newcommand{\sumstar}{\sideset{}{^*}\sum}
\newcommand{\sumprime}{\sideset{}{'}\sum}
\newcommand{\sumprimeprime}{\sideset{}{''}\sum}
\newcommand{\shortmod}{\ensuremath{\negthickspace \negthickspace \negthickspace \pmod}}
\newcommand{\V}{V\left(\frac{nm}{q^2}\right)}
\newcommand{\sumi}{\mathop{{\sum}^{\dagger}}}
\newcommand{\mz}{\ensuremath{\mathbb Z}}
\newcommand{\leg}[2]{\left(\frac{#1}{#2}\right)}
\newcommand{\muK}{\mu_{\omega}}

\title[Mean values of divisors twisted by quadratic characters]{Mean values of divisors twisted by quadratic characters}

\date{\today}
\author{Peng Gao}
\address{Department of Mathematics, School of Mathematics and Systems Science, Beihang University, P. R. China}
\email{penggao@buaa.edu.cn}

\begin{abstract}
In this paper, we evaluate the sum $\sum_{m,n}\leg {m}{n}d(n)$, where $\leg {m}{n}$ is the Kronecker symbol and $d(n)$ is the divisor
function.
\end{abstract}

\maketitle


\section{Introduction}

  As estimations for character sums have wide applications in analytic number theory, many important results are obtained in this direction. Among the many
  characters studied in the literature, quadratic Dirichlet characters receive more attentions due to their relations to ranks of elliptic curves,
class numbers, etc.  A special case of the well-known P\'olya-Vinogradov inequality (see \cite[Chap. 23]{Da}) asserts that for any positive $m \neq
\square$ (where $\square$ represents the square of a rational integer), and any $Y>0$,
\begin{align}
\label{PVquad}
    \sum_{n \leq Y} \Big (\frac {m}{n} \Big )  \ll m^{1/2}\log m,
\end{align}
   where $\leg {m}{\cdot}$ is the Kronecker symbol.

    One may regard the above P\'olya-Vinogradov inequality as a first moment estimation for quadratic Dirichlet characters. In view of this, it is natural
    to ask for similar estimations for higher moments involving quadratic character sums. For the second moment, we note the following mean square estimate
    for quadratic Dirichlet characters due to M. V. Armon \cite[Theorem 1]{Armon}:
\begin{align}
\label{squareJacobi}
   \sum_{\substack {|D| \leq X \\ D \in \mathcal{D}}} \left| \sum_{n \leq Y} \Big (\frac {D}{n} \Big ) \right|^2 \ll XY \log X,
\end{align}
   where $\mathcal{D}$ is the set of non-square quadratic discriminants. Weaker estimations were obtained earlier by M. Jutila \cites{Jutila, Jutila1} to
   study problems related to the mean values of class numbers of quadratic imaginary fields and the second moment of Dirichlet $L$-functions with primitive
   quadratic characters.

   As the P\'olya-Vinogradov inequality demonstrates certain square root cancellation, it is not possible to obtain an asymptotic expression for the
   character sum given in the left-hand side of \eqref{PVquad}. However, things change if we consider an extra average over $m \leq X$ as well. In fact, it
   is relatively easy to obtain such an asymptotic formula when $X, Y$ are far apart in size, as the P\'olya-Vinogradov inequality \eqref{PVquad} itself
   allows us to have a good control on the error term. On the other hand, things become much subtler when $X$ and $Y$ are of comparable size and it is not
   until in year 2000 that J. B. Conrey, D. W. Farmer, and K. Soundararajan \cite{CFS} determined completely the asymptotic formula for all $X,Y>0$ of the
   following sum:
\begin{align}
\label{meanchar}
    \sum_{\substack{m \leq Y \\ (m,2)=1}}\sum_{\substack{n \leq Y \\ (n,2)=1}} \Big (\frac {m}{n} \Big ).
\end{align}

   The above expression motivates us to consider similar sums involving quadratic Dirichlet characters. One can certainly do so by considering mean values
   of
any arithmetic function twisted by quadratic characters. We shall, however, in this paper, consider the following more concrete sum:
\begin{align*}
  S(X,Y)=\sum_{\substack {m \leq X \\ (m, 2)=1}}\sum_{\substack {n \leq Y \\ (n, 2)=1}} \leg {m}{n}d(n).
\end{align*}

   One reason that leads to our consideration of the above sum is that, when we square out the sums in \eqref{squareJacobi} and interchange the order of
   summations, we find that
\begin{align*}
   \sum_{\substack {|D| \leq X \\ D \in \mathcal{D}}} \left| \sum_{n \leq Y} \Big (\frac {D}{n} \Big ) \right|^2 =\sum_{n \leq Y}\sum_{\substack {|D| \leq
   X \\ D \in \mathcal{D}}} \leg{D}{n}d(n)+\sum_{\substack{n > Y \\ n=_1n_2\\ n_1 \leq Y, n_2 \leq Y}}\sum_{\substack {|D| \leq X \\ D \in \mathcal{D}}}
   \leg{D}{n}.
\end{align*}

  When one compares the first sum on the right-hand side expression above with $S(X,Y)$, it is easy to see that they differ only by certain restrictions on
  the sums and this is what gets us into the consideration of $S(X,Y)$.

  Another reason for us to consider $S(X,Y)$ comes from the work of K. Soundararajan in \cite{sound1}, in which he evaluated mollified first and second
  moments of quadratic Dirichlet $L$-functions to show that at least $87.5\%$ such $L$-functions are non-vanishing at the central point. An important part
  of the arguments in \cite{sound1} involves with the evaluation of smoothed sums of the following type:
\begin{equation*}
   \sum_{(d,2)=1}\mu^2(d) \sum_{n}\leg {8d}{n} \frac {d_j(n)}{\sqrt{n}}\omega_j\left(\frac {n}{d^{j/2}} \right )F \left(\frac {d}{X} \right), \quad j=1,2.
\end{equation*}
   Here one may regard $\omega_j, j=1,2$ and $F$ as Schwarz class functions and $d_1(n)=1, d_2(n)=d(n)$ for all $n$. Thus, other than the factor $1/\sqrt{n}$,
   the above sums when $j=1$ can be regarded as a smoothed version of the sums considered in \eqref{meanchar}, except that in the inner sum above, we
   are essentially summing over $n$ going up to the size of $d$, while in \eqref{meanchar}, the lengths of the summations are independent. In the same
   fashion, $S(X,Y)$ can be regarded as a study of the above sums for $j=2$ by making the lengths of the sums independent.

    We are now ready to state our result in the following
\begin{theorem}
\label{quadraticmean}
   For large $X$ and  $Y$, we have
\begin{align}
\label{S}
\begin{split}
   S(X,Y)=& \frac {XY^{1/2}\log^2 Y}{16}\eta (1)+XY^{1/2}P_1(\log Y) +X^{3/2}(\log X+ 2\gamma) C_1(\frac YX)+ X^{3/2}C_2(\frac YX) \\
   &+O(X^{1+\epsilon}Y^{\frac 14+\epsilon}\left ( \frac {XY^{1/2}+YX^{1/2}}{XY^{\frac 14}} \right )^{2/3}+X^{1+\epsilon}Y^{1/2+\epsilon}\left ( \frac {XY^{\frac 14}}{XY^{1/2}+YX^{1/2}} \right )^{1/3}(\frac {Y}{X})^{(1+\epsilon)/2}\\
   &+Y^{3/2+\epsilon}\left ( \frac {XY^{\frac 14}}{XY^{1/2}+YX^{1/2}} \right )^{2/3}+X^{1+\epsilon}Y^{\epsilon} \left( \frac
   {Y}{X} \right)^{2}\left ( \frac {XY^{1/2}+YX^{1/2}}{XY^{\frac 14}} \right )^{2}).
\end{split}
\end{align}
  where $\eta(s)$ is given in \eqref{eta}, $P_1(x)$ is given in \eqref{P1} and $C_1(\alpha), C_2(\alpha)$ are functions of $\alpha \geq 0$ given in \eqref{C}.
\end{theorem}

    It is easy to see that \eqref{S} gives a valid asymptotic formula when $X^{\epsilon} \ll Y \ll X^{1-\epsilon}$ for any $\epsilon>0$. We also note that via integration by parts that for $i=1,2$,
\begin{align*}
  C_i(\alpha)=O(\alpha^{3/2}\log \alpha ),  \quad \alpha \rightarrow 0.
\end{align*}
   Unlike the case studied in \cite{CFS}, the behavior of $C_i(\alpha)$ when $\alpha \rightarrow \infty$ is more complicated to analyze and we shall not worry about it here.

   Our strategy for proving Theorem \eqref{quadraticmean} is similar to that in the proof of Theorem 1 in \cite{CFS}. We first replace $S(X,Y)$ by the
   following smoothed sum
\begin{align}
\label{smoothS}
  {\mathcal S}(X,Y) &=\sum_{\substack {m \leq X \\ (m, 2)=1}}\sum_{\substack {n \leq Y \\ (n, 2)=1}} \leg {m}{n}d(n)\Phi \left(\frac {n}{Y} \right
  )W\left(\frac {m}{X} \right ).
\end{align}
    Here $\Phi$ and $W$ are smooth functions supported in $(0, 1)$, satisfying $\Phi(t)=W(t)=1$ for $t\in (\frac 1U, 1-\frac 1U)$ for a parameter $U$, and
    such that
\begin{align}
\label{UWestim}
  \Phi^{(j)}(t), W^{(j)}(t)\ll_j U^j
\end{align}
     for all integers $j \geq 0$.

   We shall apply a large sieve result for quadratic Dirichlet characters to control the size of  $S(X,Y)-{\mathcal S}(X,Y)$. We then apply the Poisson
   summation to evaluate ${\mathcal S}(X,Y)$. Choosing $U$  optimally leads to the assertion of Theorem \ref{quadraticmean}.

\section{Preliminaries}
\label{sec 2}

   In this section, we gather a few auxiliary results needed in the proof of Theorem \ref{quadraticmean}.
\subsection{Gauss sums}
\label{section:Gauss}

   For all odd
    integers $k$ and all integers $m$, we introduce the following Gauss-type
    sums as in \cite[Sect. 2.2]{sound1}
\begin{align}
\label{G}
    G_m(k)=
    \left( \frac {1-i}{2}+\left( \frac {-1}{k} \right)\frac {1+i}{2}\right)\sum_{a \shortmod{k}}\left( \frac {a}{k} \right) e \left( \frac {am}{k} \right),
\end{align}
   where $e(x)=e^{2\pi i x}$. We quote \cite[Lemma 2.3]{sound1} which determines $G_m(k)$.
\begin{lemma}
\label{lem1}
   If $(k_1,k_2)=1$ then $G_m(k_1k_2)=G_m(k_1)G_m(k_2)$. Suppose that $p^a$ is
   the largest power of $p$ dividing $m$ (put $a=\infty$ if $m=0$).
   Then for $b \geq 1$ we have
\begin{equation*}
\label{011}
    G_m(p^b)= \left\{\begin{array}{cl}
    0  & \mbox{if $b\leq a$ is odd}, \\
    \varphi(p^b) & \mbox{if $b\leq a$ is even},  \\
    -p^a  & \mbox{if $b=a+1$ is even}, \\
    (\frac {m/p^a}{p})p^a\sqrt{p}  & \mbox{if $b=a+1$ is odd}, \\
    0  & \mbox{if $b \geq a+2$}.
    \end{array}\right.
\end{equation*}
\end{lemma}

\subsection{Poisson Summation}
   For a Schwartz function $F$, we define
\begin{equation} \label{tildedef}
   \widetilde{F}(\xi)=\frac {1+i}{2}\hat{F}(\xi)+\frac
   {1-i}{2}\hat{F}(-\xi)=\int\limits^{\infty}_{-\infty}\left(\cos(2\pi \xi
   x)+\sin(2\pi \xi x) \right)F(x) \dif x,
\end{equation}
   where $\hat{F}$ denotes the Fourier transform of $F$.

    We have the following Poisson summation formula from \cite[Lemma 2.6]{sound1}:
\begin{lemma}
\label{lem2}
   Let $W$ be a smooth function with compact support on the
positive real numbers. For any odd integer $n$,
\begin{equation*}
\label{013}
  \sum_{(d,2)=1}\left( \frac {d}{n} \right)
    W\left( \frac {d}{X} \right)=\frac {X}{2n}\left( \frac {2}{n} \right)
    \sum_k(-1)^kG_k(n)\widetilde{W}\left( \frac {kX}{2n} \right),
\end{equation*}
where $\widetilde{W}$ is defined in \eqref{tildedef} and $G_k(n)$ is defined in \eqref{G}.
\end{lemma}

     We shall apply the above lemma to the function $W$ defined in the Introduction. Here we recall the following estimations on this $W$ from \cite[(4.1)]{CFS}:
\begin{align}
\label{bounds}
     \widetilde{W}^{(\mu)}(t) \ll_{j} U^{j-1}t^{-j}
\end{align}
    for all integers $\mu \geq 0, j \geq 1$ and all real $t>0$.

    On the other hand, we note that it follows from \cite[(4.2)]{CFS} that for the same $W$,
\begin{align}
\label{w0}
    \widetilde{W}(t) =\frac {1-\cos (2\pi t)+\sin (2\pi t)}{2\pi t}+O(\frac 1{U}).
\end{align}
\subsection{A large sieve for quadratic Dirichlet characters}
\label{sec2.4}
   Another important tool needed in the proof of Theorem \ref{quadraticmean} is the following large sieve inequality for quartic Dirichlet characters due
   to  D. R. Heath-Brown \cite[Theorem 1]{DRHB}.
\begin{lemma} \label{quadls}
Let $M,N$ be positive integers, and let $(a_n)_{n\in \mathbb{N}}$ be an arbitrary sequence of complex numbers. Then we have
\begin{equation*}
 \sumstar_{\substack{m \leq M}} \left| \ \sumstar_{\substack{n  \leq N}} a_n \leg{n}{m} \right|^2
 \ll_{\varepsilon} (M + N )(MN)^{\varepsilon} \sum_{n \leq N} |a_n|^2,
\end{equation*}
   for any $\varepsilon > 0$, where the asterisks indicate that $m$ and $n$ run over positive odd square-free integers.
\end{lemma}

\section{Proof of Theorem \ref{quadraticmean}}
\label{sec 3}

\subsection{Initial Reductions}

    Before we evaluate ${\mathcal S}(X,Y)$, we first want to estimate $S(X,Y)-{\mathcal S}(X,Y)$, where ${\mathcal S}(X,Y)$ is defined in \eqref{smoothS}.
    It is easy to see that $S(X,Y)-{\mathcal S}(X,Y)$ can be expressed as linear combinations of sums of the following forms:
\begin{align}
\label{Firsterr}
   \sum_{\substack {m \in I_1 \\ (m, 2)=1}}\sum_{\substack {n \in I_2 \\ (n, 2)=1}}  \leg {m}{n}d(n)H(n,m),
\end{align}
   where $H(x,y)=1$ or $\Phi \left(\frac {n}{Y} \right )W\left(\frac {m}{X} \right )$, $I_1\in \{[0, X/U], [X(1-1/U), X], [0, X] \}$, $I_2 \in \{[0, Y/U],
   [Y(1-1/U), Y], [0, Y] \}$, and the case $I_1=[0, X], I_2=[0, Y]$ is excluded.

   As the arguments are similar, we only treat the following sum:
\begin{align*}
   \sum_{\substack {m \in [X(1-1/U), X] \\ (m, 2)=1}}\sum_{\substack {n \in [0,Y] \\ (n, 2)=1}}  \leg {m}{n}d(n).
\end{align*}

   We note that
\begin{align}
\label{1stsum}
\begin{split}
  & \sum_{\substack {m \in [X(1-1/U), X] \\ (m, 2)=1}}\sum_{\substack {n \in [0,Y] \\ (n, 2)=1}}  \leg {m}{n}d(n)
   =  \sum_{\substack {a^2 \leq X \\ (a, 2)=1}}\sum_{\substack {b   \in [X(1-1/U)/a^2, X/a^2]\\ (b, 2)=1}}\sum_{\substack {e^2 \leq Y \\ (e,
   2)=1}}\sum_{\substack {f \leq Y/e^2 \\ (f, 2)=1}} \mu^2(b)\mu^2(f) \leg {a^2b}{e^2f}d(e^2f) \\
  =&  \sum_{\substack {a^2 \leq X \\ (a, 2)=1}}\sum_{\substack {e^2 \leq Y \\ (e, 2a)=1}} \sum_{\substack {b   \in [X(1-1/U)/a^2, X/a^2]\\ (b,
  2e)=1}}\sum_{\substack {f \leq Y/e^2 \\ (f, 2a)=1}} \mu^2(b)\mu^2(f) \leg {b}{f}d(e^2f).
\end{split}
\end{align}

   For fixed $a, b$, we apply the Cauchy-Schwarz inequality and Lemma \ref{quadls} to see that
\begin{align*}
  & \sum_{\substack {b   \in [X(1-1/U)/a^2, X/a^2]\\ (b, 2e)=1}}\sum_{\substack {f \leq Y/e^2 \\ (f, 2a)=1}} \mu^2(b)\mu^2(f) \leg {b}{f}d(e^2f) \\
 \ll & \left (\sum_{\substack {b   \in [X(1-1/U)/a^2, X/a^2]\\ (b, 2e)=1}}1 \right )^{1/2}\left (\sum_{\substack {b   \in [X(1-1/U)/a^2, X/a^2]\\ (b,
 2e)=1}}\left | \sum_{\substack {f \leq Y/e^2 \\ (f, 2a)=1}} \mu^2(b)\mu^2(f) \leg {b}{f}d(e^2f) \right |^2 \right )^{1/2} \\
 \ll & \sqrt{\frac {X}{Ua^2}}\left ( \left(\frac {X}{a^2}+\frac {Y}{e^2}\right )(XY)^{\epsilon}\sum_{\substack {f \leq Y/e^2 \\ (f, 2a)=1}}d^2(e^2f)
 \right )^{1/2}\ll  (XY)^{\epsilon}\sqrt{\frac {X}{Ua^2}}\left (\sqrt{\frac {XY}{a^2e^2}}+\frac {Y}{e^2}\right ).
\end{align*}
   Applying the above estimation to the right-hand side of \eqref{1stsum}, we see upon summing over $a$ and $e$ that
\begin{align*}
 \sum_{\substack {m \in [X(1-1/U), X] \\ (m, 2)=1}}\sum_{\substack {n \in [0,Y] \\ (n, 2)=1}}  \leg {m}{n}d(n)
 \ll (XY)^{\epsilon}\frac {XY^{1/2}+YX^{1/2}}{\sqrt{U}}.
\end{align*}

   One checks that the above estimation applies to all other types of sums given in \eqref{Firsterr}. We then conclude that
\begin{align}
\label{nonsmoothestm}
 S(X,Y)-{\mathcal S}(X,Y) \ll (XY)^{\epsilon}\frac {XY^{1/2}+YX^{1/2}}{\sqrt{U}}.
\end{align}

    Next, we treat ${\mathcal S}(X,Y)$ by applying Lemma \ref{lem2}, we see that
\begin{align*}
  {\mathcal S}(X,Y) &=X\sum_{n} \frac {d(n)}{2n}\leg {2}{n}\Phi \left(\frac {n}{Y} \right )
    \sum_k(-1)^kG_k(n)\widetilde{W}\left( \frac {kX}{2n} \right) \\
   & =\frac {X\widetilde{W}(0)}{2}\sum_{n } \leg {2}{n} \frac {G_0(n)d(n)}{n} \Phi \left( \frac {n}{Y} \right)+\frac {X}{2}\sum_{\substack {k \neq
   0}}(-1)^{k} \sum_{n} \leg {2}{n} \frac {G_k(n)d(n)}{n} \Phi \left( \frac {n}{Y} \right)\widetilde{W}\left( \frac {kX}{2n} \right) \\
   &:= M_{0}+M'.
\end{align*}

\subsection{The Term $M_{0}$}

    We estimate $M_{0}$ first. It follows straight from the definition that $G_0(n)=\varphi(n)$ if $n=\square$ and $G_0(n)=0$ otherwise. Thus
\begin{align*}
  M_{0}= \frac {X\widetilde{W}(0)}{2}\sum_{\substack {(n,2)=1 \\ n=\square}}\frac {\varphi(n)d(n)}{n}\Phi \left( \frac {n}{Y} \right).
\end{align*}

    By Mellin inversion, we have
\begin{align*}
    \Phi(\frac {n}{Y}) = \frac 1{2\pi i}\int\limits_{(2)} \left( \frac{Y}{n} \right)^s\widehat{\Phi}(s) \dif s,
\end{align*}
    where
\begin{align}
\label{Phihat}
  \widehat{\Phi}(s)=\int\limits^{\infty}_{0}\Phi(t)t^{s-1} \dif t.
\end{align}

     Integration by parts and using \eqref{UWestim} shows that $\widehat{\Phi}(s)$ is a function satisfying the bound for all $\Re(s) > 0$, and integers $D
     \geq 1$,
\begin{align}
\label{boundsforphi}
  \widehat{\Phi}(s) \ll (1+|s|)^{-D} U^{D-1}.
\end{align}

   We then deduce that
\begin{align}
\label{M0}
  M_{0}= \frac {X\widetilde{W}(0)}{2}\frac 1{2\pi i}\int\limits_{(2)}Y^s\widehat{\Phi}(s)\left( \sum_{\substack {(n,2)=1 \\ n=\square}}\frac
  {\varphi(n)d(n)}{n^{1+s}} \right) \dif s = \frac {X\widetilde{W}(0)}{2}\frac 1{2\pi i}\int\limits_{(2)}Y^s\widehat{\Phi}(s)\left( \sum_{\substack
  {(n,2)=1}}\frac {\varphi(n^2)d(n^2)}{n^{2+2s}} \right) \dif s .
\end{align}

   We deduce by comparing the Euler factors that
\begin{align}
\label{eta}
   \sum_{\substack {(n,2)=1 \\ n=\square}}\frac {\varphi(n)d(n)}{n^{1+s/2}}=
   \sum_{\substack {(n,2)=1 }}\frac {\varphi(n^2)d(n^2)}{n^{2+s}}=\zeta^3(s)\eta(s),
\end{align}
   where $\eta(s)=\prod_p\eta_p(s)$ with
\begin{align*}
   \eta_2(s)=(1-2^{-s})^3
\end{align*}
   and for $p>2$,
\begin{align*}
   \eta_p(s) =1-\frac 3{p^{1+s}}+  (3-4(1-\frac 1p))\frac 1{p^{2s}}-\frac 1{p^{1+3s}} .
\end{align*}
   From this we see that $\eta(s)$ is absolutely convergent in $\Re(s)>1/2$.

   We then derive from \eqref{M0} that
\begin{align*}
  M_{0}= \frac {X\widetilde{W}(0)}{2}\frac 1{2\pi i}\int\limits_{(2)}Y^s\widehat{\Phi}(s)\zeta^3(2s)\eta(2s) \dif s .
\end{align*}

    We now shift the line of integration in the above expression to the $\Re(s) =1/4+\epsilon $ line. We encounter a pole of order $3$ at $s = 1/2$. Thus
    we may write $M_0 = M_{0,1}+R_0$, where
\[  M_{0,1} = \frac {X\widetilde{W}(0)}{2}\text{Res}_{s=1/2} Y^s\widehat{\Phi}(s)\zeta^3(2s)\eta(2s) \;\; \mbox{and} \;\;
   R_1 = \frac {X\widetilde{W}(0)}{2} \frac 1{2\pi i}\int\limits_{(\frac 14+\varepsilon)}Y^s\widehat{\Phi}(s)\zeta^3(2s)\eta(2s) \dif s.
\]

     To estimate $R_1$, we note that (see \cite[Exercise 3, p. 100]{iwakow}) when $\Re(s) = 1/4+\epsilon$,
\begin{align*}
  \zeta(2s) \ll (1+|s|)^{1/4+\epsilon}.
\end{align*}

      Using this to together with \eqref{boundsforphi} by taking $D=2$, we deduce that
\begin{align}
\label{R1}
   R_1 \ll XY^{\frac 14+\epsilon}U.
\end{align}

    To evaluate $M_{0,1}$, we note the Laurent series expansion for $\zeta(s)$ at $s=1$ (see \cite[Corollary 1.16]{MVa1})
\begin{align}
\label{zetaat1}
  \zeta(s) = & \frac 1{s-1}+\gamma+O(s-1),
\end{align}
   where $\gamma$ is Euler's constant.

   We deduce from \eqref{zetaat1} the following Laurent series expansions at $s=1/2$,
\begin{align*}
  \zeta^3(2s) = & \frac 1{(2s-1)^3}+\frac {3\gamma}{(2s-1)^2}+\frac {c_1}{(2s-1)}+\ldots, \\
  \eta(2s)=& \eta(1)+2\eta'(1)(s-1/2)+2\eta''(1)(s-1/2)^2+\ldots, \\
   Y^s=& Y^{1/2}+ Y^{1/2}\log Y (s-1/2)+\frac { Y^{1/2}\log^2 Y}{2}(s-1/2)^2+\ldots, \\
   \widehat{\Phi}(s) =& \widehat{\Phi}(1/2)+\widehat{\Phi}'(1/2)(s-1/2)+\frac {\widehat{\Phi}''(1/2)}{2}(s-1/2)^2+\ldots,
\end{align*}
    where $c_1$ is an absolute constant.

   We further deduce from \eqref{Phihat} that for some absolute constants $d_1, d_2$,
\begin{align*}
   \widehat{\Phi}(1/2)=2+O(\frac 1U), \quad \widehat{\Phi}'(1/2)=d_1+O(\frac 1U), \quad \widehat{\Phi}''(1/2)=d_2+O(\frac 1U).
\end{align*}

    It follows that the residue in $M_{0,1}$ can be written as
\begin{align}
\label{P1}
   \frac {Y^{1/2}\log^2 Y}{8}\eta (1)+Y^{1/2}P_1(\log Y) +O(\frac {Y^{1/2}\log^2 Y}U),
\end{align}
  where $P_1$ is a polynomial of degree $1$, whose coefficients involve absolute constants $c_1, d_1, d_2$.

  We therefore deduce from \eqref{w0} that
\begin{align*}
   M_{0,1}=\frac {XY^{1/2}\log^2 Y}{16}\eta (1)+XY^{1/2}P_1(\log Y) +O(\frac {XY^{1/2}\log^2 Y}U).
\end{align*}

   Combining the above estimation with \eqref{R1}, we conclude that
\begin{align}
\label{M0asymp}
   M_{0}=\frac {XY^{1/2}\log^2 Y}{16}\eta (1)+XY^{1/2}P_1(\log Y) +O(\frac {XY^{1/2}\log^2 Y}U+XY^{\frac 14+\epsilon}U).
\end{align}

\subsection{The Term $M'$}
    Now suppose $k \neq 0$. By Mellin inversion, we have
\begin{align}
\label{phiW}
    \Phi \left( \frac {n}{Y} \right)\widetilde{W}\left( \frac {kX}{2n} \right) = \frac 1{2\pi i}\int\limits_{(2)} \left( \frac{Y}{n}
    \right)^s\tilde{f}(s,k) \dif s,
\end{align}
    where
\begin{align}
\label{quadraticf}
  \tilde{f}(s,k)=\int\limits^{\infty}_{0}\Phi(t)\widetilde{W}\left(\frac {kX}{2Yt}\right) t^{s-1} \dif t.
\end{align}

     Integration by parts and using \eqref{bounds} shows that $\tilde{f}(s)$ is a function satisfying the bound for all $\Re(s) > 0$, and integers $D, E>0$,
\begin{align}
\label{boundsforf}
  \tilde{f}(s,k) \ll (1+|s|)^{-D} \left( 1+\frac {|k|X}{Y} \right)^{-E+D} U^{E-1}.
\end{align}

    We apply \eqref{phiW} to recast $M'$ as
\begin{align*}
  M'=\frac {X}{2}\sum_{\substack {k \neq 0}}(-1)^{k}\int\limits_{(2)}\tilde{f}(s,k)Y^sG(1+s,k) \dif s.
\end{align*}
   where
 \begin{align*}
   G(1+s,k)=\sum_{n} \leg {2}{n} \frac {G_k(n)d(n)}{n^{1+s}}=\sum_{(n,2)=1} \frac {G_{2k}(n)d(n)}{n^{1+s}},
\end{align*}
   where the last equality above follows from the observation that, by changing variables in \eqref{G}, we have $\leg {2}{n}G_k(n) =G_{2k}(n)$ since $n$ is odd.

    Write $2k = k_1 k^2_2$ where $k_1$ is square-free
(possibly $k_1 = 1$ is the trivial character), and $k_2$ is positive. In the region
$\Re(s) > 1$, we recast $G(1+s, k)$ as
\begin{align*}
  G(1+s,k) =L \left( \frac 12+s, \chi_{k_1} \right)^2 \prod_{p}G_{p}(s,k):=L \left( \frac 12+s, \chi_{k_1} \right)^2 {\mathcal G}(s,k),
\end{align*}
   where we denote $\chi_{k}$ for the Kronecker symbol$\leg {k}{\cdot}$ and $G_{p}(s,k)$ is defined as follows:
\begin{align*}
   G_{p}(s,k)=& \left( 1-\frac {1}{p^{1/2+s}}\leg {k_1}{p} \right)^2,  \quad \text{ if $p=2$}, \\
   G_{p}(s,k)=&  \left( 1-\frac {1}{p^{1/2+s}}\leg {k_1}{p} \right)^2 \sum^{\infty}_{r=0}\frac {d(p^r)}{p^{r(1/2+s)}}\frac {G_{2k}(p^r)}{p^{r/2}},  \quad
   \text{ if $p \neq 2$}.
\end{align*}

   By Lemma \ref{lem1}, we see that for a generic $p \nmid 2k$, $G_{p}(s,k)=1-\frac 3{p^{1+2s}}+\frac {2}{p^{3/2+3s}}\leg {k_1}{p}$, this shows that
   ${\mathcal G}(s,k)$ is holomorphic in
 $\Re(s) > 0$. From our evaluation of $G_{p}(s,k)$ for $p \nmid 2k$ we see that for  $\Re{s} \geq \epsilon$,
\begin{align*}
   {\mathcal G}(s,k) \ll |k|^{\epsilon}\prod_{p | 2k}|G_{p}(s,k)|.
\end{align*}

   To derive a bound for $G_{p}(s,k)$ when $p | 2k$, we suppose that $p^a \| 2k$. By the trivial bound $|G_k(p^r)| \leq p^r$, it follows that $|G_{p}(s,k)| \leq (a+1)^2$ for these $p$. We then conclude that for  $\Re{s}
   \geq \epsilon$,
\begin{align}
\label{boundsforG}
   {\mathcal G}(s,k) \ll |k|^{\epsilon}.
\end{align}

   Using this, we see that
\begin{align*}
   M' =\frac {X}{2}\sum_{\substack {k \neq 0}}(-1)^{k}\int\limits_{(2)}\tilde{f}(s,k)Y^sL \left( \frac 12+s, \chi_{k_1} \right)^2 {\mathcal G}(s,k) \dif s
   .
\end{align*}

     We now move the line of integration to the line $\Re(s) = \epsilon$ to see that we encounter poles only when $k_1 = 1$ (so that $L(s, \chi_{k_1}) =
     \zeta(s)$) at $s=1/2$, in which case there is a pole of order $2$ at $s = 1/2$. Thus we may write $M'_2 = M_{1,2}+R_2$, where (by an obvious change of
     notation, writing $2k^2$ in place of the corresponding $k$ and observing that $k^2_1=k^2_2$ if and only if $k_1=\pm k_2$)
\[   M_{1,2} = X\displaystyle{\text{Res}_{s=1/2}} \sum^{\infty}_{\substack {k =1}} Y^{s}\tilde{f}\left( \frac1{2}, 2k^2 \right){\mathcal G}
\left( \frac 1{2}, 2k^2 \right) \zeta \left( \frac 12+s \right)^2, \]
   and
\begin{align}
\label{R2}
R_2 = \frac {X}{2}\sum_{\substack {k \neq 0}}(-1)^{k} \frac 1{2\pi i}\int\limits_{(\varepsilon)}\tilde{f}(s,k)Y^s L \left( \frac 12+s, \chi_{k_1}
\right)^2 {\mathcal G}(s,k) \dif s.
\end{align}

     To evaluate $M_{1,2}$, we note that it follows from \eqref{zetaat1} that near $s=1/2$,
\begin{align*}
  \zeta(1/2+s)^2 = & \frac 1{(s-1/2)^2}+\frac {2\gamma}{s-1/2}+O(s-1/2).
\end{align*}

    It follows that
\begin{align}
\label{M12}
\begin{split}
  M_{1,2} = &  XY^{1/2} \log Y\sum^{\infty}_{\substack {k =1}} \tilde{f}\left( \frac1{2}, 2k^2 \right){\mathcal G} \left( \frac 1{2}, 2k^2
  \right) \\
  & + XY^{1/2}\sum^{\infty}_{\substack {k =1}} \tilde{f}\left( \frac1{2}, 2k^2 \right){\mathcal G}\left( \frac 1{2}, 2k^2 \right)\left (2\gamma+
  \frac {\tilde{f}'\left( \frac1{2}, 2k^2 \right)}{\tilde{f}\left( \frac1{2}, 2k^2 \right)}+ \frac {{\mathcal G}' \left( \frac 1{2}, 2k^2
  \right)}{{\mathcal G} \left( \frac 1{2}, 2k^2 \right)}\right ).
\end{split}
\end{align}

   We further note that when $s=1/2$, we have by Lemma \ref{lem1} that
\begin{align*}
   \sum^{\infty}_{r=0}\frac {d(p^r)}{p^{r(1/2+s)}}\frac {G_{4k^2}(p^r)}{p^{r/2}} >1.
\end{align*}
   We deduce from this easily that
\begin{align*}
   \frac {{\mathcal G}' \left( \frac 1{2}, 2k^2 \right)}{{\mathcal G} \left( \frac 1{2}, 2k^2 \right)} \ll |k|^{\epsilon}.
\end{align*}
   We apply the definition of $\tilde{f}\left( s, k \right)$ given in \eqref{quadraticf} together with the bound \eqref{boundsforG} to see that
\begin{align}
\label{sumfg}
\begin{split}
   \sum^{\infty}_{\substack {k =1}} \tilde{f}\left( \frac1{2}, 2k^2 \right){\mathcal G} \left( \frac 1{2}, 2k^2 \right) =& \int\limits^{\infty}_{0}\frac
   {\Phi(t)}{\sqrt{t}}\sum^{\infty}_{\substack {k =1}}\widetilde{W}\left(\frac {k^2X}{Yt}\right){\mathcal G} \left( \frac 1{2}, 2k^2 \right)  \dif t, \\
   \sum^{\infty}_{\substack {k =1}} \tilde{f}'\left( \frac1{2}, 2k^2 \right){\mathcal G} \left( \frac 1{2}, 2k^2 \right) =& \int\limits^{\infty}_{0}\frac
   {\Phi(t)\log t}{\sqrt{t}}\sum^{\infty}_{\substack {k =1}}\widetilde{W}\left(\frac {k^2X}{Yt}\right){\mathcal G} \left( \frac 1{2}, 2k^2 \right)  \dif t,
   \\
   \sum^{\infty}_{\substack {k =1}} \tilde{f}\left( \frac1{2}, 2k^2 \right){\mathcal G}' \left( \frac 1{2}, 2k^2 \right) =& \int\limits^{\infty}_{0}\frac
   {\Phi(t)}{\sqrt{t}}\sum^{\infty}_{\substack {k =1}}\widetilde{W}\left(\frac {k^2X}{Yt}\right){\mathcal G}' \left( \frac 1{2}, 2k^2 \right)  \dif t.
\end{split}
\end{align}

   We then apply \eqref{bounds} with $j=1$ for $|k| \geq \sqrt{UY/X}$ and \eqref{w0} for smaller $k$ to see that
\begin{align}
\label{WG}
\begin{split}
  \sum^{\infty}_{\substack {k =1}}\widetilde{W}\left(\frac {k^2X}{Yt}\right){\mathcal G} \left( \frac 1{2}, 2k^2 \right)
  =& \frac {Yt}{2\pi X}  \sum^{\infty}_{\substack {k =1}} \frac {{\mathcal G} \left( \frac 1{2}, 2k^2 \right)}{k^2}\left ( 1-\cos (\frac {2\pi k^2
  X}{Yt})+\sin (\frac {2\pi k^2 X}{Yt}) \right )+ O(U^{-(1-\epsilon)/2}(\frac {Y}{X})^{(1+\epsilon)/2}), \\
  \sum^{\infty}_{\substack {k =1}}\widetilde{W}\left(\frac {k^2X}{Yt}\right){\mathcal G}' \left( \frac 1{2}, 2k^2 \right)
  = &\frac {Yt}{2\pi X}  \sum^{\infty}_{\substack {k =1}} \frac {{\mathcal G}' \left( \frac 1{2}, 2k^2 \right)}{k^2}\left ( 1-\cos (\frac {2\pi k^2
  X}{Yt})+\sin (\frac {2\pi k^2 X}{Yt}) \right )+ O(U^{-(1-\epsilon)/2}(\frac {Y}{X})^{(1+\epsilon)/2}).
\end{split}
\end{align}

    Applying \eqref{sumfg} and \eqref{WG} to \eqref{M12}, we deduce that
\begin{align}
\label{M12asymp}
  M_{1,2} =  X^{3/2}(\log X+ 2\gamma) C_1(\frac YX)+ X^{3/2}C_2(\frac YX)+O(XY^{1/2}\log YU^{-(1-\epsilon)/2}(\frac {Y}{X})^{(1+\epsilon)/2}+\frac
  {Y^{3/2}\log Y}{U}),
\end{align}
    where we define for $\alpha \geq 0$,
\begin{align}
\label{C}
\begin{split}
   C_1(\alpha)
  = & \frac {1}{2\pi}\sum^{\infty}_{\substack {k =1}} \frac {{\mathcal G} \left( \frac 1{2}, 2k^2 \right)}{k^2}\int\limits^{\alpha}_{0}\sqrt{y}  \left (
  1-\cos (\frac {2\pi k^2}{y})+\sin (\frac {2\pi k^2 }{y}) \right )dy, \\
  C_2(\alpha)
  = & \frac {1}{2\pi}\sum^{\infty}_{\substack {k =1}} \frac {{\mathcal G} \left( \frac 1{2}, 2k^2 \right)}{k^2}\int\limits^{\alpha}_{0}\log y \sqrt{y}
  \left ( 1-\cos (\frac {2\pi k^2}{y})+\sin (\frac {2\pi k^2 }{y}) \right )dy \\
  & +\frac {1}{2\pi}\sum^{\infty}_{\substack {k =1}} \frac {{\mathcal G}' \left( \frac 1{2}, 2k^2 \right)}{k^2}\int\limits^{\alpha}_{0}\sqrt{y}  \left (
  1-\cos (\frac {2\pi k^2}{y})+\sin (\frac {2\pi k^2 }{y}) \right )dy .
\end{split}
\end{align}

     To estimate $R_2$, we note that (see \cite[Exercise 3, p. 100]{iwakow}) when $\Re(s) =\epsilon$,
\begin{align*}
  L \left( \frac 12+s, \chi_{k_1} \right) \ll (1+|s|)^{1/4+\epsilon}.
\end{align*}

      Using this together with \eqref{boundsforf} by taking $D=2, E=4$ for all $k$ and \eqref{boundsforG}, we deduce from \eqref{R2} that
\begin{align}
\label{R2bound}
   R_2 \ll XY^{\epsilon}  \left( \frac {Y}{X} \right)^{2}U^{3} \sum_{|k|\geq 1} \frac {|k|^{\varepsilon}}{|k|^{2}}  \ll   XY^{\epsilon} \left( \frac
   {Y}{X} \right)^{2}U^{3}.
\end{align}

   We then conclude from \eqref{M12asymp} and \eqref{R2bound} that
\begin{align}
\label{M'asymp}
  M'=X^{3/2}(\log X+ 2\gamma) C_1(\frac YX)+ X^{3/2}C_2(\frac YX)+O(XY^{1/2}\log YU^{-(1-\epsilon)/2}(\frac {Y}{X})^{(1+\epsilon)/2}+\frac
  {Y^{3/2}\log Y}{U}+XY^{\epsilon} \left( \frac
   {Y}{X} \right)^{2}U^{3}.
\end{align}

\subsection{Conclusion}
  Combining \eqref{nonsmoothestm}, \eqref{M0asymp} and \eqref{M'asymp}, we see that

\begin{align*}
   S(X,Y)=& \frac {XY^{1/2}\log^2 Y}{16}\eta (1)+XY^{1/2}P_1(\log Y) +X^{3/2}(\log X+ 2\gamma) C_1(\frac YX)+ X^{3/2}C_2(\frac YX) \\
   &+O(\frac {XY^{1/2}\log^2 Y}U+XY^{\frac 14+\epsilon}U+(XY)^{\epsilon}\frac {XY^{1/2}+YX^{1/2}}{\sqrt{U}}+XY^{1/2}\log YU^{-(1-\epsilon)/2}(\frac {Y}{X})^{(1+\epsilon)/2} \\
   &+\frac
  {Y^{3/2}\log Y}{U}+XY^{\epsilon} \left( \frac
   {Y}{X} \right)^{2}U^{3}).
\end{align*}

    By taking
\begin{align*}
    U=\left ( \frac {XY^{1/2}+YX^{1/2}}{XY^{\frac 14}} \right )^{2/3},
\end{align*}
   we derive the expression \eqref{S} in Theorem \ref{quadraticmean} and this completes the proof.

\noindent{\bf Acknowledgments.} The author is supported by NSFC grant 11871082.

\bibliography{biblio}

\begin{bibdiv}
\begin{biblist}

\bib{Armon}{article}{
      author={Armon, M.~V.},
       title={Averages of real character sums},
        date={1999},
     journal={J. Number Theory},
      volume={77},
      number={2},
       pages={209\ndash 226},
}

\bib{CFS}{article}{
      author={Conrey, J.~B.},
      author={Farmer, D.~W.},
      author={Soundararajan, K.},
       title={Transition mean values of real characters},
        date={2000},
     journal={J. Number Theory},
      volume={82},
      number={1},
       pages={109\ndash 120},
}

\bib{Da}{book}{
      author={Davenport, H.},
       title={{M}ultiplicative {N}umber {T}heory},
     edition={Third edition},
      series={Graduate Texts in Mathematics},
   publisher={Springer-Verlag},
     address={Berlin, etc.},
        date={2000},
      volume={74},
}

\bib{DRHB}{article}{
      author={Heath-Brown, D.~R.},
       title={A mean value estimate for real character sums},
        date={1995},
     journal={Acta Arith.},
      volume={72},
      number={3},
       pages={235\ndash 275},
}

\bib{iwakow}{book}{
      author={Iwaniec, H.},
      author={Kowalski, E.},
       title={{A}nalytic {N}umber {T}heory},
      series={American Mathematical Society Colloquium Publications},
   publisher={American Mathematical Society},
     address={Providence},
        date={2004},
      volume={53},
}

\bib{Jutila1}{article}{
      author={Jutila, M.},
       title={On character sums and class numbers},
        date={1973},
     journal={J. Number Theory},
      volume={5},
       pages={203\ndash 214},
}

\bib{Jutila}{article}{
      author={Jutila, M.},
       title={On the mean value of ${L}(1/2,\chi)$ for real characters},
        date={1981},
     journal={Analysis},
      volume={1},
      number={2},
       pages={149\ndash 161},
}

\bib{MVa1}{book}{
      author={Montgomery, H.~L.},
      author={Vaughan, R.~C.},
       title={{M}ultiplicative number theory. {I}. {C}lassical theory},
      series={Cambridge Studies in Advanced Mathematics},
   publisher={Cambridge University Press},
     address={Cambridge},
        date={2007},
      volume={97},
}

\bib{sound1}{article}{
      author={Soundararajan, K.},
       title={Nonvanishing of quadratic {D}irichlet {$L$}-functions at
  $s=\frac{1}{2}$},
        date={2000},
     journal={Ann. of Math. (2)},
      volume={152},
      number={2},
       pages={447\ndash 488},
}

\end{biblist}
\end{bibdiv}
\bibliographystyle{amsxport}



\end{document}